\documentclass[11pt]{article}
\usepackage{amsmath,amssymb,amsthm}

\newcommand{\Mbar}{\overline{\mathcal{M}}}

\newtheorem{definition}{Definition}
\newtheorem{lemma}[definition]{Lemma}
\newtheorem{theorem}[definition]{Theorem}
\newtheorem{corollary}[definition]{Corollary}

\newtheorem*{theorem'}{Theorem}

\begin{document}

\title{Grafting hyperbolic metrics and Eisenstein series}         
\author{Kunio Obitsu \and Scott A. Wolpert\footnote{Primary: 32G15; Secondary: 32Q15, 30F60, 53D30}}
\date{November 8, 2007}       
\maketitle

\begin{abstract}
The family hyperbolic metric for the plumbing variety $\{zw=t\}$ and the non holomorphic Eisenstein series $E(\zeta;2)$ are combined to provide an explicit expansion for the hyperbolic metrics for degenerating families of Riemann surfaces. Applications include an asymptotic expansion for the Weil-Petersson metric and a local form of symplectic reduction. 
\end{abstract}

\section{Introduction}

A genus $g$, $d$ punctured Riemann surface with $2g-2+d$ positive admits a complete hyperbolic structure.  A family of marked hyperbolic structures degenerates provided there are closed geodesics with lengths tending to zero.  Degenerating families provide an important means for investigating global quantities and the global geometry of the moduli space of Riemann surfaces.  Degenerating families also provide a means for investigating families of structures with boundary lengths tending to zero. Hyperbolic structures with geodesic boundary can be doubled across boundaries to form structures without boundary.  The examinations are based on expansions for degenerating structures.  

L. Bers considered certain degenerating families of marked hyperbolic structures \cite{Bersdeg}.  He found that away from collars about short geodesics the structures limit to hyperbolic structures with cusps with neighborhoods of cusps removed.   Expansions have been presented for such degenerating families of hyperbolic structures.  We combine the approaches of 
\cite{Wfinf,Wlhyp} to present an explicit expansion in terms of the hyperbolic metrics (see (\ref{hyp}) and (\ref{hyp2}) below) of fibers of the plumbing family $\mathcal V_0$ and the non holomorphic Eisenstein series 
$E(\zeta;2)$.  Beginning with a Riemann surface $R$ with a formal pairing of $n$ pairs of punctures (additional unpaired punctures of $R$ are permitted), we use the fiberwise gluing of $\mathcal V_0$ to $R-\{cusp\ neighborhoods\}$ to construct a degenerating family 
$\{R_t\}$.  We then use geometric-interpolation to combine the $\mathcal V_0$ fiber metrics with the $R$ metric to form grafted metrics $dg_t^2$.  A counterpart interpolation the {\em melding} is introduced for the Eisenstein series.  The construction differs from Thurston's grafting of projective structures for Riemann surfaces \cite{Gdproj}.  We then use uniqueness properties for the prescribed curvature equation to show that the grafted metric and meldings combine to provide an expansion.

\begin{theorem'}
For a grafting cutoff function $\eta$ and all $t$ small the hyperbolic metric $ds_{hyp}^2$ of $R_t$ has the expansion
\[
ds_{hyp}^2=dg_t^2\bigl(1+\frac{4\pi^4}{3}\sum\limits_{k=1}^n (\log |t_k|)^{-2}(E^{\dagger}_{k,1}+E^{\dagger}_{k,2})+O(\sum\limits_{k=1}^n(\log|t_k|)^{-3})\bigr).
\] 
The functions $E_{k,1}^{\dagger}$ and $E_{k,2}^{\dagger}$ are the meldings of the Eisenstein series $E(\cdot;2)$ associated to the pair of cusps plumbed to form the $k^{th}$ collar.  The $O$-term refers to the intrinsic $C^1$-norm of a  function on $R_t$. The bound depends on the choice of $\eta$ and a lower bound for the injectivity radius for the complement of the cusp regions in $R$.
\end{theorem'}
A length expansion is a first application.  For the plumbing family $\{R_t\}$ there is a comparison between fibers by domain inclusion of collar complements.  The comparison provides a setting for relating lengths for geodesics disjoint from the collars
\[
\ell_{\alpha}(\ell)=\ell_{\alpha}(0)+\sum\limits_{k=1}^n \frac{\ell^2_k}{6}\int_{\alpha}(E_{k,1}+E_{k,2})ds+
O(\sum\limits_{k=1}^n \ell_k^3).
\]
Similarly in \cite{Wldis} Wolf's infinite energy harmonic map degeneration expansion was employed with Eisenstein series to study the variation of eigenvalues of the Laplacian.  

A further application is for the K\"{a}hler and symplectic geometry of the Teichm\"{u}ller space $\mathcal T$ of  punctured Riemann surfaces with hyperbolic metrics 
\cite{Ahsome,Busbook,ImTan,Ngbook}.  The infinitesimal deformations of a surface $R$ are represented by the Beltrami differentials $\mathcal H (R)$ harmonic with respect to the hyperbolic metric \cite{Ahsome}. The cotangent space of $\mathcal T$ at $R$ is $Q(R)$ the space of holomorphic quadratic differentials with at most simple poles at the punctures of $R$.  The holomorphic tangent-cotangent pairing is
\[
(\mu,\varphi)=\int_R \mu \varphi
\]
for $\mu\in\mathcal H(R)$ and $\varphi\in Q(R)$.  Elements of $\mathcal H (R)$ are symmetric tensors given as $\overline\varphi(ds^2)^{-1}$.  A Hermitian form for $\mathcal H(R)$ defines a Hermitian metric for Teichm\"{u}ller space.  For $\mu,\nu\in\mathcal H(R)$ and the hyperbolic area element $dA$ the Weil-Petersson (WP) form is
\[
\langle\mu,\nu\rangle_{WP}=\int_R\mu\bar\nu\, dA, \qquad \cite{Ahsome}
\]
and for $R$ with punctures $p_1,\dots,p_d$ and corresponding Eisenstein series $E_1,\dots,E_d$ the Takhtajan-Zograf (TZ) form for the puncture $p_j$ is
\[
\langle\mu,\nu\rangle_{TZ,p_j}=\int_R\mu\bar\nu E_j\,dA, \qquad \cite{TZpunc}.
\]
The Takhtajan-Zograf metric is $\sum\limits^d_{j=1}\langle\ ,\ \rangle_{TZ,p_j}$.

The WP metric is K\"{a}hler and non complete with non pinched negative curvature.  The metric has been studied extensively; see \cite{Wlcomp,Wlpers} for beginning references.  The WP metric provides a complete $CAT(0)$ geometry for the augmented Teichm\"{u}ller space $\overline{\mathcal T}$ \cite{DW2,Wlcomp,Yam2}. The mapping class group is the WP isometry group and the augmented Teichm\"{u}ller space is the closed convex hull of the set of totally noded Riemann surfaces \cite{MW,Wlcomp}.  Degenerating families of hyperbolic metrics have been used to develop an asymptotic expansion for the metric and its curvature 
\cite{DW2,Zh2,Msext,LSY1,LSY2,Schu,Wlhyp,Wlcomp}.  

The TZ metric is K\"{a}hler and non complete \cite{Obit1,Obit2,TZpunc}.  L. Takhtajan and P. Zograf introduced the metric in their investigation of the local index formula for $\overline{\partial}$ and the Quillen metric \cite{TZpunc,TZ}. The authors showed for punctured Riemann surfaces that the Chern forms for the determinant line bundles for families of $\overline{\partial}$-operators are sums of the WP and TZ K\"{a}hler forms \cite{TZ}.  K. Obitsu, W.-K. To and L. Weng have used degenerating families of hyperbolic metrics to develop an asymptotic expansion for the TZ metric \cite{OTW}.   In the Deligne-Mumford moduli space of stable curves transverse to the divisor $\mathcal D$ of noded curves the TZ metric is smaller than the WP metric. 

We apply the above hyperbolic metric expansion to refine the asymptotic expansion for the WP metric.  Prior expansions \cite{Msext,OTW,Wlcomp} focus on directions transverse to the divisor $\mathcal D$.  We consider directions parallel to $\mathcal D$.  For the family $\{R_t\}$ we express the expansion in terms of geodesic length with the expansion $\ell=2\pi^2(-\log|t|)^{-1}+O((\log|t|)^{-2})$ \cite[Ex. 4.3]{Wlhyp}. 

\begin{theorem'}
In the Deligne-Mumford moduli space of stable curves for tangent subspaces parallel to the divisor of noded curves the WP metric has the expansion for all small $\ell$
\[
\langle\ ,\ \rangle_{WP}^{parallel}(\ell)=\\
\langle\ ,\ \rangle_{WP}(0) + \sum_{k=1}^n\frac{\ell_k^2}{3}
\bigl(\langle\ ,\ \rangle_{TZ,k,1} + \langle\ ,\ \rangle_{TZ,k,2}\bigr) + O(\sum\limits_{k=1}^n\ell_k^{3}\bigr)
\]
for the TZ forms for punctures paired to form the $k^{th}$ node.
\end{theorem'}

There are considerations for WP geometry.  The augmentation locus of $\overline{\mathcal T}$ is a union of closed convex strata \cite{Wlcomp,Yam2}.  A complete convex subset of a $CAT(0)$ space is the base of a projection \cite[Chap. II.2]{BH}.  A stratum in $\overline{\mathcal T}$ is described as the $0$-level set for a collection $\sigma$ of geodesic lengths $\{\ell_k\}_{k\in\sigma}$.  Distance to the stratum is given as $(2\pi\sum_{k\in\sigma}\ell_k)^{1/2}+O(\sum_{k\in\sigma}\ell_k^{5/2})$ \cite[Coro. 21]{Wlcomp} and \cite[Coro. 4.10]{Wlbhv} with the projection distance non increasing.  Accordingly the perturbation in $\sum_{k\in\sigma}\ell_k$ of the WP metric from the $0$-level set is expected to be non negative.  From the above parallel to the level set the initial perturbation is strictly positive.   

There are considerations for K\"{a}hler forms based on a second characterization of TZ forms.  For a holomorphic family of punctured Riemann surfaces the family of cotangent spaces along a puncture $p$ defines a tautological holomorphic line bundle $\psi_p$ over the base of the family.  The hyperbolic metric can be renormalized to define a metric $\|\ \|_{ren,p}$ for the cotangent space at the puncture.  The Chern form for the family of cotangent spaces along the puncture $p$ is given in terms of the TZ exterior $2$-form $c_1(\|\ \|_{ren,p})=\frac43\omega_{TZ,p}$ \cite{Wlcusps}. The asymptotic expansion for K\"{a}hler forms becomes
\begin{equation}
\label{kexp}
\omega(\ ,\ )_{WP}^{parallel}(\ell)=\\
\omega(\ ,\ )_{WP}(0) + \sum_{k=1}^n\frac{\ell_k^2}{4}
\bigl(c_1(\psi_{k,1}) + c_1(\psi_{k,2})\bigr) + O(\sum\limits_{k=1}^n\ell_k^{3}\bigr).
\end{equation}

An application is for the volume of moduli spaces of bordered Riemann surfaces.  M. Mirzakhani considered the moduli space $\mathcal M(b_1,\dots,b_d)$ of  bordered Riemann surfaces with geodesic boundary components of prescribed lengths $(b_1,\dots, b_d)$ \cite{Mirvol}. The spaces $\mathcal M(b_1,\dots,b_d)$ form an $(\mathbb R_{\ge 0})^d$ bundle over $\mathcal M$ the moduli space of genus $g$, $d$ punctured Riemann surfaces.  A Riemann surface with geodesic boundary and a {\em point} on each boundary is alternately described by a $d$ punctured Riemann surface and a product of $d$ factors of $S^1$, a principal torus-bundle over a punctured Riemann surface. In 
\cite[Secs. 3 \& 4]{Mirwitt} Mirzakhani used  symplectic reduction (for an $(S^1)^d$ quasi-free action following Guillemin-Sternberg) to provide a simple description for the family of WP K\"{a}hler forms
\[
\omega_{WP}\big\vert_{\mathcal M (b_1,\dots,b_d)}=\omega_{WP}\big\vert_{\mathcal M }+ \sum\limits_j \frac{b_j^2}{4}\,c_1(\psi_j).
\]
The expansion is a relation for cohomology classes on $\overline{\mathcal M}$ and is presented for a general choice of principal connection for the $S^1$ bundles.  Formula (\ref{kexp}) is a local form of the symplectic reduction.  

\section{The plumbing family hyperbolic metric}
For a Riemann surface $R$ with a finite area hyperbolic metric and punctures there is a natural 
{\em cusp coordinate} (with  germ unique modulo rotation) at each puncture: at a puncture $p$, the coordinate $z$ with $z(p)=0$ and the hyperbolic metric of $R$ given as $(\frac{|dz|}{|z|\log |z|})^2$.  For the surface uniformized by the upper half-plane $\mathbb H$ with $p$ represented by a width-one cusp at infinity then a loop once encircling the puncture corresponds to the deck transformation $\zeta \rightarrow\zeta +1$; the natural cusp coordinate is given as $z=e^{2\pi i \zeta}$ for $\zeta$ the natural $\mathbb H$ coordinate.  

We recall the construction of the plumbing family  \cite[esp. Sec. 2.4]{Wlhyp} and the grafting of hyperbolic metrics \cite[esp. Sec. 3.3]{Wlhyp}.  Considerations begin with the 
{\em plumbing variety} 
$\mathcal{V}=\{(z,w,t)\mid zw=t,\ |z|, |w|, |t|<1\}$.  
The defining function $zw-t$ has differential $z\,dw+w\,dz-dt$.  Consequences
are that $\mathcal{V}$ is a smooth variety, $(z,w)$ are global coordinates,
while $(z,t)$ and $(w,t)$ are not.  Consider the projection
$\Pi:\mathcal{V}\rightarrow D$ onto the $t$-unit disc.  The projection $\Pi$ is a submersion,
except at $(z,w)=(0,0)$;  we consider $\Pi:\mathcal{V}\rightarrow D$ as a
(degenerate) family of open Riemann surfaces.  The $t$-fiber, $t\ne 0$, is the
hyperbola germ $zw=t$ or equivalently the annulus
$\{|t|<|z|<1,\,w=t/z\}=\{|t|<|w|<1,\,z=t/w\}$.  The $0$-fiber is the
intersection of the unit ball with the union of the coordinate axes in
$\mathbb{C}^2$; on removing the origin the union becomes
$\{0<|z|<1\}\cup\{0<|w|<1\}$.  Each fiber of 
$\mathcal{V}_0=\mathcal{V}-\{0\}\rightarrow D$ has a complete hyperbolic metric:

\begin{equation}
\begin{split}
\label{hyp}
\mbox{for }t\ne0,&\mbox{ on }\{|t|<|z|<1\} \mbox{ then}\\
&ds^2_t=\Bigl(\frac{\pi}{\log |t|} \csc \frac{\pi\log |z|}{\log |t|}\Bigl|\frac{dz}{z}\Bigr|\Bigr)^2;\\
\mbox{for }t=0,&\mbox{ on }\{0<|z|<1\}\cup \{0<|w|<1\} \mbox{ then}\\
&ds^2_0=\Bigl(\frac{|d\zeta|}{|\zeta|\log|\zeta|}\Bigr)^2\mbox{ for } \zeta=z,\,w.
\end{split}
\end{equation}
The family of hyperbolic metrics $(ds_t^2)$ is a continuous metric, degenerate
only at the origin, for the vertical line bundle of $\mathcal{V}$.  In
particular we have the elementary expansion

\begin{equation}
\begin{split}
\label{hyp2}
ds^2_t=&\, \Bigl(\frac{|d\zeta|}{|\zeta|\log|\zeta|}\Bigr)^2\,\Bigl(\Theta\csc\Theta\Bigr)^2\mbox{\quad for\quad}\Theta=\frac{\pi\log |z|}{\log |t|}\\
=&\, ds_0^2\,\bigl(1+\frac13\Theta^2+\frac{1}{15}\Theta^4+\dots\bigr).
\end{split}
\end{equation}

Consider $R$ a finite union of Riemann surfaces with punctures. A plumbing family is the fiberwise gluing of the complement of cusp neighborhoods in $R$ and the plumbing variety $\mathcal V=\{(z,w,t)\mid zw=t, |z|, |w|, |t|<1\}$.  For a positive constant $c_*<1$ and initial surface $R$, with puncture $p$ with cusp coordinate $z$ and puncture $q$ with cusp coordinate $w$, we construct a family $\{R_t\}$.  For $|t|<c_*^4$ the resulting surface $R_t$ will be independent of $c_*$; the constant $c_*$ will serve to specify the overlap of coordinate charts and to define a {\em collar} in each $R_t$.

We first describe the gluing of fibers.  For $|t|<c_*^4$, remove from $R$ the punctured discs $\{0<|z|\le |t|/c_*\}$ about $p$ and $\{0<|w|\le |t|/c_*\}$ about $q$ to obtain a surface $R^*_{t/c_*}$ \cite[Sec. 2.4, second parag.]{Wlhyp}.  For $t\ne 0$, form an identification space $R_t$, by identifying the annulus $\{|t|/c_*<|z|<c_*\}\subset R^*_{t/c_*}$ with the annulus $\{|t|/c_*<|w|<c_*\}\subset R^*_{t/c_*}$ by the rule $zw=t$.  The resulting surface $R_t$ is the {\em plumbing} for the prescribed value of $t$.  We note for $|t|<|t'|$ that there is an inclusion of $R^*_{t/c_*}$ in $R^*_{t'/c_*}$; the inclusion maps provide a way to compare structures on the surfaces.  The inclusion maps are a basic feature of the plumbing construction. We next describe the plumbing family.  Consider the variety $\mathcal V_{c_*}=\{(z,w,t)\mid zw=t, |z|,|w| < c_*, |t|<c_*^4\}$ and the disc $D_{c_*}=\{|t|<c_*^4\}$.  
The complex manifolds $M=R^*_{t/c_*}\times D_{c_*}$ and $\mathcal V_{c_*}$ have holomorphic projections to the disc $D_{c_*}$.  The variables $z,w$ denote prescribed coordinates on $R^*_{t/c_*}$ and on $\mathcal V_{c_*}$.  There are holomorphic maps of subsets of $M$ to $\mathcal V_{c_*}$, commuting with the projections to $D_{c_*}$, as follows
\[
(z,t)\stackrel{\hat F}{\rightarrow}(z,t/z,t) \mbox{ and } (w,t)\stackrel{\hat G}{\rightarrow}(w,t/w,t).
\]
The identification space $\mathcal F=M\cup\mathcal V_{c_*}/\{\hat F,\hat G \mbox{ equivalence}\}$ is the {\em plumbing family} $\{R_t\}$ with projection to $D_{c_*}$ (an analytic fiber space of Riemann surfaces in the sense of Kodaira \cite{Kod}).  For $0<|t|<c_*^4$, the $t$-fiber of $\mathcal F$ is the surface $R_t$ constructed by overlapping annuli.

The hyperbolic metric $ds^2$ of $R$ restricted to $R^*_{c_*}$ can be combined with the fiber hyperbolic metric $ds_t^2$ 
of $\mathcal V$ restricted to $\mathcal V_{c_*^2}$ to provide an approximation,  the grafted metric, to the hyperbolic of $R_t$, 
\cite[Sec. 3.3]{Wlhyp}.  To this purpose select $\eta(a)$, a smooth function, identically unity for $a\le a_0< 0$ and identically zero 
for $a\ge 0$. (We further restrict $t$ to satisfy $e^{2a_0}c_*^2\ge |t|$.) For the variable $z$ denoting the cusp coordinate of $R^*_{t/c_*}$ and the fiber coordinate for $\mathcal V_{c_*}$ define a smooth metric in a neighborhood of $\{|z|=c_*\}$ by the rule
\[
dg_t^2=(ds^2)^{1-\eta}(ds_t^2)^{\eta} \mbox{ with } \eta=\eta(\log (|z|/c_*));
\]
and similarly in a neighborhood of $\{|w|=c_*\}$ (for a Riemann surface a metric is a section of an $\mathbb R^+$-bundle with structure group $\mathbb R^{\times}$; the geometric-interpolation is well defined).  In the {\em collar bands} $\{e^{a_0}c_*\le|z|\le c_*\}$ and $\{e^{a_0}c_*\le|w|\le c_*\}$ the {\em grafted metric} $dg_t^2$ is the geometric-interpolation of the two metrics.  Away from the collar the grafted metric is prescribed as $ds^2$, the R hyperbolic metric, and in the complement of the bands in the collar the grafted metric is prescribed as $ds_t^2$.   

The  grafted metric is the combination of the metric $ds^2$ for $R^*_{t/c_*}$ and the fiber metric $ds_t^2$ for $\mathcal V_{c_*^2}$.  The basic relation $ds_t^2=ds_0^2(\Theta\csc \Theta)^2$ for $\Theta=\frac{\pi \log |z|}{\log |t|}$ for the fiber metrics of $\mathcal V$ provides that $ds_t^2$ is real analytic in $\epsilon=\frac{\pi}{\log |t|}$, locally for nonzero $z$.  For $\zeta$ the coordinate for $\mathbb H$ and $\Theta=\frac{-2\pi^2\Im \zeta}{\log |t|}$ the relation provides $ds^2_t=ds_0^2(1+\frac{4\pi^4}{3}(\log |t|)^{-2}(\Im \zeta)^2+O((\log |t|)^{-4}))$ a presage of the expansion in the next section.  The construction of the grafted metric is motivated by two important  properties: for $t$ small, $dg_t^2=ds^2(1+O((\log |t|)^{-2}))$, and in the collar core $dg_t^2$ coincides with the fiber metrics for $\mathcal V$. Importantly the curvature of the grafted metric is identically $-1$ on the complement of the collar bands, and satisfies in the collar bands
\begin{equation}
\label{graftcurv}
K_{\mbox{{graft}}}=-1-\frac{\epsilon^2}{6}\Lambda\ + \ O(\epsilon^4)
\end{equation}
for $\Lambda=\frac{\partial}{\partial a}(a^4\frac{\partial}{\partial a} \eta)$, for $a=\log |z|,\,\log |w|$, \cite[Sec. 3.4.MG, Def'n 3.8]{Wlhyp}.  The construction of the grafted metric is local for a pair of punctures.  The considerations apply for plumbing multiple pairs of punctures.  

\section{The special truncation of $E(\zeta;2)$}

The Eisenstein series $E(\zeta;s)$ for $s=2$ plays a central role in our considerations.  We begin with standard estimates.  For a Fuchsian group $\Gamma$, with the stabilizer of infinity the group of integer-translations $\Gamma_{\mathbb Z}$, consider the relative Poincar\'{e} series
\[
E(\zeta;2)=\sum\limits_{A\in \Gamma_{\mathbb Z}\backslash \Gamma}(\Im A(\zeta))^2, \cite{Brl,Vn}.
\]
The function $(\Im \zeta)^2$ on $\mathbb H$ is an eigenfunction of the hyperbolic Laplacian with eigenvalue 2. The above series is associated to the cusp at infinity.  In general eigenfunctions $\varphi$ on $\mathbb H$ with eigenvalue $\lambda$ satisfy a mean-value property; consequently for $\epsilon_0$ positive an eigenfunction $\varphi$ satisfies
\[
|\varphi(\zeta_0)|\le c(\epsilon_0,\lambda)\int_{B(\zeta_0;\epsilon_0)}|\varphi|d\mathcal A
\]
for $B$ the metric ball and $d\mathcal A$ the hyperbolic area element \cite{Fay}.  A  consequence is that the Eisenstein series converges locally uniformly on $\mathbb H$.  Uniform bounds for the Eisenstein series are also a further consequence.  To this purpose we consider the injectivity radius of the projection to $\mathbb H/\Gamma$ of a point $\zeta$.  The cusp region $\{\Im \zeta>1\}$ is  {\em precisely invariant} by $\Gamma$: the translate of the region by an element of $\Gamma$ is either disjoint from the region, or coincides with the region and the transformation is an integer-translation.  Accordingly for $\Im \zeta >1$ the injectivity radius of the image point on $\mathbb H/\Gamma$ is approximately $(2\Im \zeta)^{-1}$.  Accordingly a point of $\mathbb H/\Gamma$ with injectivity radius at least $\epsilon_0$ has lifts  with $\Gamma$-orbit contained in $\{2\Im \zeta\le \epsilon_0^{-1}\}$; the bound is uniform among all Fuchsian groups with stabilizer of infinity $\Gamma_{\mathbb Z}$. Accordingly there is a universal bound at points with injectivity radius at least $\epsilon_0$
\begin{multline*}
E(\zeta_0;2)\le c(\epsilon_0)\int_{B(\zeta_0;\epsilon_0/2)}Ed\mathcal A\\
\le c(\epsilon_0)\int_{\cup_{A\in \Gamma_{\mathbb Z}\backslash\Gamma}A(B(\zeta_0;\epsilon_0/2))}(\Im \zeta)^2d\mathcal A
\le c(\epsilon_0)\int_{2\Im \zeta\le e^{\epsilon_0/2}\epsilon_0^{-1};\,0\le \Re \zeta \le 1}d\mathcal E
\end{multline*}
for $d\mathcal E$ the Euclidean area element.  We are interested in bounds for all cusp regions of $\mathbb H/\Gamma$  uniform in $\Gamma$.  By successive conjugations of $\Gamma$, the cusps of $\Gamma$ can be represented  at infinity with stabilizer $\Gamma_{\mathbb Z}$.  For a cusp so represented at infinity, the quotient space $\{\Im \zeta>1\}/\Gamma_{\mathbb Z}$ embeds in $\mathbb H/\Gamma$ to determine a {\em cusp region}.  Cusp regions for distinct cusps are disjoint.  

The Eisenstein series, transformed for a cusp represented at infinity, has the expansion
\[
E(\zeta;2)=\delta_{\infty}(\Im \zeta)^2+\hat e(\zeta)
\]
with $\delta_{\infty}$ the indicator for the defining cusp for the series, and $\hat e(\zeta)$ bounded as $O((\Im \zeta)^{-1})$ for $\Im \zeta$ large \cite{Brl,Vn}.  We write $E^{\star}$ for the discontinuous function that satisfies $E^{\star}(\zeta;2)=E(\zeta;2)-(\Im \zeta)^2$ for the cusp region defining the series represented at infinity and that otherwise agrees with $E$.

\begin{lemma}
\label{Ebound}
For $\zeta$ the $\mathbb H$ coordinate and a cusp represented by a width-one cusp at infinity the derivatives of order at most $k$ of $E^{\star}$ for $\Im \zeta > 2$ are bounded as $O((\Im \zeta)^{-1})$.  The  bound depends on $k$ but is independent of $\Gamma$.
\end{lemma}
\begin{proof} The supremum bound is the main consideration.  For the metric ball of radius 
$\frac12$, the mean-value estimate provides that $E^{\star}$ is suitably bounded on $\Im \zeta =\frac32$.  Accordingly for a suitable positive constant $c$ we have that $-c(\Im \zeta)^{-1}\le E^{\star}\le c(\Im \zeta)^{-1}$ on $\Im \zeta=\frac32$ with $E^{\star}$ and $(\Im \zeta)^{-1}$ each vanishing at infinity.  The pair of functions $E^{\star}$ and $(\Im \zeta)^{-1}$ are eigenfunctions of the hyperbolic Laplacian with the positive eigenvalue 2.  By the maximum principle $E^{\star}-c(\Im \zeta)^{-1}$, resp. $E^{\star}+c(\Im \zeta)^{-1}$, realizes its maximum, resp. minimum, on the boundary of the domain $\{\Im \zeta\ge\frac32\}$, \cite{GT}.  Realized at infinity the maximum , resp. minimum, is zero.  The supremum bound $-c(\Im \zeta)^{-1}\le E^{\star}\le c(\Im \zeta)^{-1}$ for $\Im \zeta \ge \frac32$ is a consequence.  The bound for the derivatives of $E^{\star}$ will be a consequence of the supremum-bound interior  estimates \cite[Chap. 5]{BJS} and the eigenequation.  In particular for $B\subset\mathbb R^2$, a bounded domain, and $B_0\subset\subset B$ there exists a constant $c=c(B,B_0)$ such that for $\Delta$ the Euclidean Laplacian then $\|P\|_{1,B_0}\le c(\|P\|_{0,B}+\|\Delta P\|_{0,B})$, where $\|\ \|_k$ is the Euclidean $C^k$-norm for a domain, \cite[Chap. 5]{BJS}.  To make the desired analysis, Euclidean translate $E^{\star}$ to a fixed ball.  The $C^1$-bound immediately follows from the eigenequation.  More generally, for $\Delta_m=\frac{\partial^m}{\partial \zeta^p\partial\bar\zeta^{m-p}}$, a derivative of order $m$, we have that $\Delta\Delta_mE^{\star}=\Delta_m\Delta E^{\star}=\Delta_m2(\Im \zeta)^{-2}E^{\star}$.  It follows that $\Delta\Delta_mE^{\star}$ is bounded in terms of the derivatives of $E^{\star}$ of order at most $m$.  We proceed by induction and apply the interior Schauder estimate to obtain bounds for the derivatives of order at most $m+1$ in the cusp.  \end{proof}

We are now ready to introduce a special truncation of the Eisenstein series given a specification of $\eta$ and the parameters $c_*$ and $t$ as above.   

\begin{definition}
\label{trun}
The special truncation $E^\#$ of the Eisenstein series is given by modification in the cusp regions.  
For the cusp defining the series where 
\[
E(\zeta;2)=(\Im \zeta)^2+\hat e(\zeta),
\]
define for $\Im \zeta >1$ and $\chi=1-\eta$,
\[
E^{\#}(\zeta;2)=\chi(-2\pi \Im \zeta-\log c_*)(\Im \zeta)^2+\chi(-2\pi\Im\zeta+(\log c_*/|t|)+a_0)\hat e(\zeta),
\]
and for a remaining cusp represented at infinity define for $\Im \zeta >1$, 
\[
E^{\#}(\zeta;2)=\chi(-2\pi\Im\zeta+\log (c_*/|t|)+a_0)E(\zeta;2).
\]
\end{definition}
It is important that for the change of variables $-2\pi\Im\zeta=\log|z|$, the function $\frac{\partial}{\partial a}\chi(-2\pi \Im \zeta-\log c_*)$ has support in the band $\{e^{a_0}c_*<|z|<c_*\}$  and the function $\frac{\partial}{\partial a}\chi(-2\pi\Im\zeta+\log (c_*/|t|)+a_0)$ has support in the band $\{|t|/c_*<|z|<e^{-a_0}|t|/c_*\}$.

We return to the plumbing construction. Consider $R$ a finite union of Riemann surfaces with punctures and $E$ the Eisenstein series for a cusp $p$.  Consider the family $\{R_t\}$ given by plumbing $n$ pairs of cusps.  Each surface $R_t$, $t=(t_1,\dots,t_n)$, is obtained by identifying $n$ pairs of subannuli for cusp regions. We first extend the special truncation $E^\#$ by zero on the components of $R$ not containing $p$.   We are ready to define the {\em melding} $E^{\dagger}$ of Eisenstein series and analyze $(D_{{graft}}-2)E^{\dagger}$.   For cusp coordinates $z,w$ of $R$ and constant $c_*<1$ the punctured discs $\{0<|z|\le|t|/c_*\}$, $\{0<|w|\le|t|/c_*\}$ are removed and  the annuli $\{|t|/c_*<|z|<c_*\}$, $\{|t|/c_*<|w|<c_*\}$ are identified by the rule $zw=t$ to form a collar.  
For $z=e^{2\pi i \zeta},\,\zeta\in\mathbb H,$ the identified annulus is covered by $\{\log(|t|/ c_*)<-2\pi\Im \zeta< \log c_*\}$; the {\em primary collar band} $\{c_*e^{a_0}\le |z|\le c_*\}$ is covered by the strip $\{\log c_*+a_0\le -2\pi \Im \zeta\le \log c_*\}$; the {\em secondary collar band} $\{|t|/c_*\le |z|\le |t|e^{-a_0}/c_*\}$ is covered by the strip $\{\log (|t|/c_*)\le -2\pi \Im \zeta\le \log (|t|/c_*)-a_0\}$.   
The extended $E^\#$ has support in the $z,w$ cusp regions contained in $\{|z|\ge|t|/c_*\}\cup\{|w|\ge|t|/c_*\}$.  We  define the {\em melding} $E^{\dagger}$ on $R_t$ of the Eisenstein series associated to $p$ to be: given on each pair of identified annuli by the sum of values of $E^\#$ at $z$ and at $w=z/t$ on the overlap of $\{|t|/c_*<|z|<c_*\}$ and $\{|t|/c_*<|w|<c_*\}$, and equal to $E^\#$ on the complement of the identified annuli.  The melding is a smooth function.

\begin{lemma}
\label{Emeld}
Relative to the grafted metric and the prescribed coordinates the melding $E^{\dagger}$ satisfies
\[
(D_{{graft}}-2)E^{\dagger}=\frac{-1}{4\pi^2}\Lambda+O((-\log|t|)^{-1})
\] 
with the remainder term referring to the $C^0$-norm of a function supported in the collar. 
\end{lemma}
\begin{proof} The grafted metric is given as $dg_t^2=(ds^2)^{1-\eta}(ds_t^2)^{\eta}$ for $ds^2$ the hyperbolic metric of $R$ and $ds_t^2$ the fiber hyperbolic metric of $\mathcal V$.  We analyze the contributions of the metrics $dg_t^2$, $ds_t^2$ and $ds^2$ and the contributions of the truncation of the Eisenstein series.  On the collar complement the grafted metric coincides with $ds^2$ and $E^{\dagger}=E$ is annihilated by the operator. The quantity $(D_{{graft}}-2)E^{\dagger}$ is supported in the collar.   

We consider separately the contribution of $E^\#$ on the $z$ and $w$ annuli. 
We first evaluate $(D_{{graft}}-2)E^\#$ on the cusp region defining the series.   On the primary and secondary collar bands we have uniformly that $dg_t^2=ds_t^2(1+O((\log|t|)^{-2}))$; $D_{{graft}}E^\#$ can be replaced on the collar bands with $D_tE^\#$ and a  remainder.  The remainder for the primary band is suitably bounded, and the expression for the secondary band will be considered below. We find with the relations $ds_t^2=ds_0^2(\Theta\csc\Theta)^2$ for $\Theta=\frac{\pi\log|z|}{\log |t|},\, z=e^{2\pi i\zeta},\, \zeta =x+iy,$ the hyperbolic Laplacian $D=y^2(\frac{\partial^2}{\partial x^2}+\frac{\partial^2}{\partial y^2})$, and the expansion from Definition \ref{trun} for the truncated Eisenstein series for $\Im \zeta>1$ 
\[
E^\#=(1-\eta_1)y^2+(1-\eta_2)\hat e(\zeta)
\]
that
\begin{gather*}
(D_t-2)E^\#=(\Theta^{-1}\sin \Theta)^2DE^\#-2E^\#=\\
 (\Theta^{-1}\sin \Theta)^2\bigl(2E^\#-4y^3\frac{\partial \eta_1}{\partial y}-y^4\frac{\partial^2 \eta_1}{\partial y^2} -2y^2\frac{\partial \hat e(\zeta)}{\partial y}\frac{\partial \eta_2}{\partial y}-y^2\hat e(\zeta)\frac{\partial^2 \eta_2}{\partial y^2}\bigr)-2E^\#\\
=A+B+C
\end{gather*}
for 
\begin{gather*}
A=((\Theta^{-1}\sin\Theta)^2-1)2E^\#,\quad B=-(\Theta^{-1}\sin \Theta)^2(4y^3\frac{\partial \eta_1}{\partial y}+y^4\frac{\partial^2 \eta_1}{\partial y^2})
\end{gather*}
and
\begin{gather*}
C=-(\Theta^{-1}\sin \Theta)^2(2y^2\frac{\partial \hat e(\zeta)}{\partial y}\frac{\partial \eta_2}{\partial y}+y^2\hat e(\zeta)\frac{\partial^2 \eta_2}{\partial y^2}\bigr).
\end{gather*}
We consider the quantities $A$, $B$ and $C$ in order.  From Lemma \ref{Ebound} $E^\#$ has magnitude $O(y^{-1})$ for $y$ large, and since $-2\pi y=\log |z|$, $\Theta=\frac{\pi\log|z|}{\log|t|}$ then $E^\#$ is bounded 
as $O((\log|t|)^{-1}\Theta^{-1})$.  Since $\Theta^{-1}((\Theta^{-1}\sin\Theta)^2-1)$ is bounded for $0<\Theta<\pi$, it follows that $A$ has magnitude $O((-\log|t|)^{-1})$, as desired.  The support of $B$ is contained in the primary 
collar band where $(\Theta^{-1}\sin\Theta)^2=1+O((\log|t|)^{-2})$. Thus the quantity $B$ is given as 
$-\frac{\partial}{\partial y}(y^4\frac{\partial \eta_1}{\partial y})$ and a suitable remainder.  
Since $-2\pi y=\log|z|$ 
it further follows from (\ref{graftcurv}) that $B$ is given as $\frac{-1}{4\pi^2}\Lambda$ and a suitable remainder, as desired.  
The support of the third term $C$ is contained in the secondary collar band since the derivatives 
of $\eta_2$ are factors.  On the support, the difference between $-2\pi y=\log|z|$ and $\log|t|$ is bounded; on the support  $\Theta=\pi+O((-\log|t|)^{-1})$, and $(\Theta^{-1}\sin\Theta)^2$ is bounded as $O((\log|t|)^{-2})$. By 
Lemma \ref{Ebound} on the support of a derivative of $\eta_2$ the quantity $(2y^2 \frac{\partial \hat e}{\partial y}+y^2\hat e)$ is bounded as $O(-\log|t|)$.  The quantity $C$ is bounded as $O((-\log|t|)^{-1})$, as desired.   

On a remaining cusp region for identification $E^\#=\hat e(\zeta)$ with $\hat e(\zeta)$ bounded as $O(y^{-1})$ and in the above notation $(D_t-2)E^\#=A+C$.  The above considerations apply to provide that the evaluation is bounded as $O((-\log|t|)^{-1})$.  The melding $E^{\dagger}$ is given as the sum of values of $E^\#$ for pairs of cusp regions and thus has the desired expansion. \end{proof}

We continue to consider $R$ a finite union of Riemann surfaces with punctures and $\{R_t\}$ the family given by plumbing $n$ pairs of cusps.  For $t=(t_1,\dots,t_n),\,|t_k|<c^4_*$, the Riemann surface $R_t$ is independent of the  parameter $c_*$.  We will use the intrinsic $C^1$-norm: for $f$ a differentiable function set $\|f\|_{{C^1}}=\sup_R|f|+\sup_{\mathbf T R}|(ds^2)^{-1/2}df|$.  We present the degeneration expansion for the hyperbolic metric.

\begin{theorem}
\label{hypexp}
Given a choice of $c_*<1$ and a cutoff function $\eta$, then for all $t$ small the hyperbolic metric $ds_{{hyp}}^2$ of $R_t$ has the expansion
\[
ds_{hyp}^2=dg_t^2\bigl(1+\frac{4\pi^4}{3}\sum\limits_{k=1}^n (\log |t_k|)^{-2}(E^{\dagger}_{k,1}+E^{\dagger}_{k,2})+O(\sum\limits_{k=1}^n(\log|t_k|)^{-3})\bigr).
\] 
The functions $E_{k,1}^{\dagger}$ and $E_{k,2}^{\dagger}$ are the meldings of the Eisenstein series $E(\cdot;2)$ associated to the pair of cusps plumbed to form the $k^{{th}}$ collar.  The $O$-term refers to the intrinsic $C^1$-norm of a  function on $R_t$. The bound depends on the choice of $c_*$, $\eta$ and a lower bound for the injectivity radius for the complement of the cusp regions in $R$.
\end{theorem}
\begin{proof} The expansion is the combination of Expansion 4.2 \cite[pg. 445]{Wlhyp} with expansion (\ref{graftcurv}) and Lemma \ref{Emeld} above.  We follow the approach of the {\em curvature correction equation} \cite[Sec. 4.2]{Wlhyp}.  The approach remains valid for surfaces with punctures and compactly supported curvature perturbations.  From equation (4.3) and Lemma 4.1 of \cite{Wlhyp} $ds_{hyp}^2=dg_t^2(1+2(D_{{graft}}-2)^{-1}(1+K_{{graft}})+O_k(\|1+K_{{graft}}\|^2))$.  
The remainder terms for (\ref{graftcurv}) and Lemma \ref{Emeld} above are supported in the plumbing collar bands; for a collection of plumbings the combined remainder terms have compact support away from the collars in the plumbed surface.  From Estimate A.2 and A.3 of 
\cite[pg.467]{Wlhyp} the $C^1$-intrinsic norm of the operator $(D_{{graft}}-2)^{-1}$ applied to the remainder terms has norm bounded by the $C^0$-norm of the remainder terms and vanishes at any additional punctures (the bound does not depend on the injectivity radius of the support).  Now combining expansion (\ref{graftcurv}) and Lemma \ref{Emeld} with the uniqueness of solutions of $(D_{{graft}}-2)u=v$ with $u$ vanishing at punctures, we find $2(D_{{graft}}-2)^{-1}(1+K_{{graft}})=\sum_k\frac{4\pi^4}{3}(\log|t_k|)^{-2}(E_{k,1}^{\dagger}+E_{k,2}^{\dagger}+O((-\log|t_k|)^{-1}))$.  
\end{proof}

{\bf Comparison to the infinite energy harmonic map expansion.}  For sake of comparison it is instructive to write the metric expansion in terms of the length $\ell=2\pi^2(-\log|t|)^{-1}+O((\log|t|)^{-2})$, \cite[Ex. 4.3]{Wlhyp}, of the geodesic of the $ds_t^2$ fiber metric of $\mathcal V$
\[
ds_{{hyp}}^2=dg_t^2\bigl(1+\sum_k\frac{\ell_k^2}{3}(E_{k,1}^{\dagger}+E_{k,2}^{\dagger})+O(\sum_k\ell_k^3)\bigr).
\]
M. Wolf's infinite energy harmonic map parameterization of the plumbing family provides that
\[
ds_{{hyp}}^2=ds^2+\sum_k\ell_k^2\bigl(\frac12\Re(\Phi_{k,1}+\Phi_{k,2})+\frac16(E_{k,1}+E_{k,2})ds^2\bigr)+O(\sum_k\ell_k^3),
\]
where $\Phi_{k,1}$ and $\Phi_{k,2}$ are holomorphic quadratic differentials with unit {\em residue} for the cusps paired for the $k^{{th}}$-node \cite[esp. Cor. 5.4]{Wfinf}.  The two expansions are  similar in spite of the different contexts.  The grafting expansion is given in terms of a metric on each surface $R_t$ and the comparison between metrics is by conformal inclusion maps.  Wolf's expansion is in terms of a family of metrics on $R$ and the comparison between fibers is by infinite energy (non conformal) harmonic maps.  The appearance of holomorphic quadratic differentials is attributed to a difference of gauge.  We found in 
\cite{Wldis} for a closed geodesic $\alpha$ of $R$ and $\Phi$ the holomorphic wt. 4 Eisenstein series (the $\Gamma_{\mathbb Z}$-coset sum of $(d\zeta)^2$ for the cusp represented at infinity) that
\[
\int_{\alpha}E(\zeta;2)ds=3\int_{\alpha}\Phi(ds)^{-1}.
\]
Each approach for the expansion gives for the length of a closed geodesic disjoint from plumbing collars
\[
\ell_{\alpha}(\ell)=\ell_{\alpha}(0)+\sum_k\frac{\ell^2_k}{6}\int_{\alpha}(E_{k,1}+E_{k,2})ds+O(\sum_k\ell_k^3),
\]
in agreement with \cite[formula (1.9)]{Wldis}. The positivity of the variation of length is a basic property of degeneration; see also expansion (\ref{hyp2}).  The positive order $2$ variation in $\ell$ corresponds to the positive remainder term for the pairings of WP gradients of length functions, \cite{Wlbhv}.

\section{The WP degeneration expansion}

We apply the expansion for the hyperbolic metric to develop a degeneration expansion for the WP metric.  We identify the leading term, a general remainder term third-order in distance, and for directions parallel to the divisor of nodal curves we identify the TZ metric as the fourth-order term in distance.  We first describe {\em local coordinates} for the deformation space of a Riemann surface with nodes and punctures. We employ a
modification of the standard coordinates \cite{Bersdeg,Msext} and combine the
formulation and the expansion for the hyperbolic metric to present a refined form of Masur's expansion of the WP metric.

We continue with discussion of the plumbing variety $\mathcal V$ and consider infinitesimal variations of the 
fibers.  To illustrate a general point we first take $c_*=1$ and consider a fiber 
$\{|t|\le|z|\le 1\}$.  
The parameter $t$ is a boundary point of the annulus.  The
boundary points $t,\,1$ will be included in the data for gluings.  To describe
the variation of annuli with boundary points, we now specify a quasiconformal
map $\zeta$ from the pointed $t$-annulus to the pointed $t'$-annulus
$\zeta(z)=zr^{\beta(r,t')},\, z=re^{i\theta},$ with
$\frac{\partial\beta}{\partial r}$ compactly supported in the annulus.  The
boundary conditions are $\zeta(1)=1$, and by specification
$t|t|^{\beta(|t|,t')}=t'$.  On differentiating in $t'$ and evaluating at
$(|t|,t)$ we find the boundary condition $t\log |t|\,\dot{\beta}(|t|,t)=1$.
The infinitesimal variation of the map is the vector field
$\dot{\zeta}(z)=z\log r\, \dot{\beta}(r,t)$ for $\dot{\zeta}\,, \dot{\beta}$
the first $t$-derivatives.   The map $\zeta$ varies from the identity and has
Beltrami differential 

\begin{equation}
\label{belt}
\bar{\partial}\dot{\zeta}=\frac{z}{2\bar{z}}\,\frac{\partial}{\partial\log r}(\dot{\beta}(r,t)\log r)\,\frac{\overline{dz}}{dz}.
\end{equation}
For sake of later application we evaluate the pairing with a quadratic
differential $z^{\alpha}\bigl(\frac{dz}{z}\bigr)^2$,
\begin{equation}
\label{pair}
\begin{split}
\int_{\{|t|<|z|<1\}}\bar{\partial}\dot{\zeta}\,z^{\alpha}\bigr(\frac{1}{z}\bigl)^2 d\mathcal E&=\int_{\{|t|<|z|<1\}}\frac{z^{\alpha}}{2z\bar{z}}\frac{\partial}{\partial \log r}(\dot{\beta}\log r)\,d\mathcal E \\
\mbox{where for }&\alpha=0 \mbox{, then}\\
&=\,\pi \dot{\beta}\log r\Big|^1_{|t|}=\frac{-\pi}{t},\\
\mbox{and other}&\mbox{wise, then}\\
&=0,\\
\end{split}
\end{equation}
for $d\mathcal E$ the Euclidean area element and where we have applied the boundary
condition for $\dot{\beta}$; the evaluation involves fixing a normalization
for the Serre duality pairing and agrees with \cite[Prop. 7.1]{Msext}.

We recall the description of  families of {\em Riemann surfaces with nodes} 
\cite{Bersdeg,Msext,Wlhyp}.  A Riemann surface $R$ with nodes and punctures is a connected
complex space, such that every point has a neighborhood isomorphic to either
the unit disc in $\mathbb{C}$, the punctured unit disc in $\mathbb C$, or the germ at the origin in $\mathbb{C}^2$ of
the union of the coordinate axes.  The points with punctured disc neighborhoods are 
{\em primary punctures.} $R$ is {\em stable} provided each component
of $R-\{nodes\}$ has negative Euler characteristic, i.e. has a hyperbolic
metric.  A regular $k$-differential on $R$ is the assignment of a meromorphic
$k$-differential $\Theta_j$ for each component $R_j$ of $R-\{nodes\}$ such
that: i) each $\Theta_*$ has poles only at the punctures of $R_*$ with orders
at most $k$, ii) poles at primary punctures have order at most $k-1$, 
and iii) if punctures $p,\,q$ are paired to form a node then
$Res_p\Theta_*=(-1)^k\,Res_q\Theta_*$ \cite{Bersdeg}.

We recall the deformation theory of Riemann surfaces with punctures and then
with nodes. A deformation neighborhood of the marked surface
$R$ is given by specifying smooth Beltrami differentials $\nu_1,\dots,\nu_m$
spanning the Dolbeault group $H^{0,1}_{\bar{\partial}}(\overline
R,\mathcal{E}((\kappa p_1\cdots p_d)^{-1}))$ for $\kappa$ the canonical bundle
of $\overline{R}$ (the closure of $R$) and $p_1,\dots,p_d$ the point line bundles for the punctures 
\cite{Kod}.  For $s\in \mathbb{C}^m$ set $\nu(s) = \sum_ks_k\nu_k$; for $s$
small there is a Riemann surface $R^{\nu(s)}$ and a diffeomorphism
$\zeta:R\rightarrow R^{\nu(s)}$ satisfying $\bar{\partial}\zeta=\nu(s)
\partial\zeta$.  The family of surfaces $\{R^{\nu(s)}\}$ represents a
neighborhood of the marked Riemann surface in its Teichm\"{u}ller space.  
The Beltrami differentials can be modified a small amount so that in terms of each {\em cusp coordinate}
the diffeomorphisms $\zeta^{\hat\nu(s)}$ are simply rotations \cite[Lemma 1.1]{Wlspeclim};
$\zeta^{\hat\nu(s)}$is a hyperbolic isometry in a neighborhood of the cusps;
$\zeta^{\hat\nu(s)}$ cannot be complex analytic in $s$, but is real analytic.
We note that for $s$ small the $s$-derivatives of $\nu(s)$ and
$\hat{\nu}(s)$ are close. We say that $\zeta^{\hat{\nu}(s)}$ {\em preserves
cusp coordinates}.   The parameterization provides a key ingredient for
obtaining simplified estimates of the degeneration of hyperbolic metrics and
an improved expansion for the WP metric.
 
We describe a {\em local manifold cover} of the compactified
moduli space $\Mbar$.  For $R$ having nodes and punctures, $R_0=R-\{nodes\}$ is a union of
Riemann surfaces with punctures.  The quasiconformal deformation space of
$R_0$, $Def(R_0)$, is the product of the Teichm\"{u}ller spaces of the
components of $R_0$.  As above for $m= \dim\,Def(R_0)$ there is a real analytic family of Beltrami differentials
$\hat{\nu}(s)$, $s$ in a neighborhood of the origin in $\mathbb{C}^m$, such
that $s\rightarrow R_s=R^{\hat{\nu}(s)}$ is a coordinate parameterization of a
neighborhood of $R_0$ in $Def(R)$ and the prescribed mappings
$\zeta^{\hat{\nu}(s)}:R_0\rightarrow R^{\hat{\nu}(s)}$ preserve the cusp
coordinates at each puncture.  For $R$ with $n$ nodes we  prescribe
the plumbing data $(U_k,V_k,z_k,w_k,t_k),\,k=1,\dots,n$, for
$R^{\hat{\nu}(s)}$, where $z_k$ on a neighborhood $U_k$ and $w_k$ on a neighborhood $V_k$ are cusp
coordinates relative to the $R^{\hat{\nu}(s)}$-hyperbolic metric (the plumbing
data varies with $s$).  The parameter $t_k$ parameterizes opening the
$k^{th}$  node.  For all $t_k$ suitably small, perform the $n$ prescribed
plumbings to obtain the family $R_{s,t}=R^{\hat{\nu}(s)}_{t_1,\dots,t_n}$.
The tuple $(s,t)=(s_1,\dots,s_m,t_1,\dots,t_n)$ provides real analytic local
coordinates, the {\em hyperbolic metric plumbing coordinates}, for the local
manifold cover of $\Mbar$ at $R$, \cite{Msext,Wlcut} and 
\cite[Secs. 2.3,2.4]{Wlhyp}. The coordinates have a special property: for $s$ fixed the
parameterization is holomorphic in $t$.  The property is a basic feature of
the plumbing construction.  The family $R_{s,t}$ parameterizes the small
deformations of the marked noded surface with punctures $R$.

We review the geometry of the local manifold covers.  For a complex manifold
$M$ the complexification $\boldsymbol{T}^{\mathbb{C}} M$ of the
$\mathbb{R}$-tangent bundle  is decomposed into the subspaces of {\em
holomorphic} and {\em antiholomorphic} tangent vectors.   A Hermitian metric
$g$ is prescribed on the holomorphic subspace.  For a general complex
parameterization  $s=u+iv$ the coordinate $\mathbb{R}$-tangents are expressed
as $\frac{\partial}{\partial u}=\frac{\partial}{\partial
s}+\frac{\partial}{\partial\bar{s}}$ and $\frac{\partial}{\partial
v}=i\frac{\partial}{\partial s}-i\frac{\partial}{\partial\bar{s}}$.   For the
$R_{s,t}$ parameterization the $s$-parameters are not holomorphic while for
$s$-parameters fixed the $t$-parameters are holomorphic;
$\{\frac{\partial}{\partial
s_j}+\frac{\partial}{\partial\bar{s}_j},\,i\frac{\partial}{\partial
s_j}-i\frac{\partial}{\partial\bar{s}_j},\,\frac{\partial}{\partial t_k},\,i
\frac{\partial}{\partial t_k}\}$ is a basis over $\mathbb{R}$ for the tangent
space of the local manifold cover.  For a smooth Riemann surface the dual of
the space of holomorphic tangents is the space of quadratic differentials with at most simple poles at punctures.
The following is a modification of Masur's result \cite[Prop. 7.1]{Msext}.

\begin{lemma}
\label{coords}
The hyperbolic metric plumbing coordinates (s,t) are real analytic and for s
fixed the parameterization is holomorphic in t. Provided the modification
$\hat{\nu}$ is small, for a neighborhood of the origin there are families in
$(s,t)$ of regular 2-differentials $\varphi_j$, $\psi_j$, j=1,\dots,m and
$\eta_k$, k=1,\dots,n  such that:   \begin{enumerate}
\item  Each regular $2$-differential has an expansion of the form $\varphi(s,t)=\varphi(s,0)+O(t)$ locally away from the nodes of $R$. 
\item For $R_{s,t}$ with $t_k\ne 0$, all k, $\{\varphi_j,\psi_j,\eta_k\,,i\eta_k\}$ forms the dual basis to  
$\{\frac{\partial\hat{\nu}(s)}{\partial s_j}+\frac{\partial\hat{\nu}(s)}{\partial \bar s_j},i\frac{\partial\hat{\nu}(s)}{\partial s_j}-i\frac{\partial\hat{\nu}(s)}{\partial \bar s_j}, \frac{\partial}{\partial t_k},i\frac{\partial}{\partial t_k}\}$ over $\mathbb{R}$.
\item For $R_{s,t}$ with $t_k=0$, all k, the $\eta_k$, k=1,\dots,n, are trivial and the   $\{\varphi_j,\psi_j\}$ span the dual of the holomorphic subspace $\boldsymbol{T}Def(R_0)$.
\end{enumerate}
\end{lemma}

\begin{proof} The setting compares to that considered by Masur.  The new element:
the variation of the plumbing data with $s$ is prescribed by a Schiffer variation for a
gluing-function real analytically depending on the parameter $s$ \cite[pg.410]{Wlcut}.  As noted  for $s$ fixed, plumbing produces a holomorphic
family.  Following Masur the families of regular $2$-differentials
$\{\varphi_j,\psi_j,\eta_k\}$ are obtained by starting with a local frame
$\mathcal{F}$ of regular $2$-differentials and prescribing the  pairings with
$\{\frac{\partial\hat{\nu}}{\partial s_j},\frac{\partial\hat{\nu}}{\partial
\bar s_j},\frac{\partial}{\partial t_k}\}$  \cite[Sec. 5 and Prop. 7.1]{Msext}. At an initial point the  basis is simply given by a linear
transformation of the frame $\mathcal{F}$.  Also the elements of $\mathcal F$ each have
an expansion of the desired form in $t$.  The prescribed basis will 
exist in a neighborhood provided the pairings are continuous with the desired expansion in $t$.  We first
consider the pairings with $\frac{\partial}{\partial t_k}$. From  (\ref{belt})
we have the Beltrami differential for the  pairing with $\frac{\partial}{\partial t_k}$,
$k=1,\dots,n$.  In particular for a plumbing collar of $R_{s,t}$ let $z$ (or
$w$) be the coordinate of the plumbing.  A quadratic differential $\varphi$ on
$R_{s,t}$ can be factored on the collar into a product of
$\bigl(\frac{dz}{z}\bigr)^2$ and a function holomorphic in $z$.  We write
$\mathcal{C}_k(\varphi)$ for the constant coefficient of the Laurent expansion
of the function factor.  From (\ref{pair}) the pairing with
$\frac{\partial}{\partial t_k}$ is the linear functional
$-\frac{\pi}{t_k}\mathcal{C}_k$.  From Masur's considerations \cite[Sec. 5, esp. 5.4, 5.5]{Msext} the pairing of $\frac{\partial}{\partial t_k}$ with the local frame $\mathcal F$ 
is continuous with the desired expansion,
and there are  regular $2$-differentials $\{\varphi_j,\psi_j,\eta_k^*\}$ with:
$\mathcal C_{\ell}(\varphi_j)=\mathcal C_{\ell}(\psi_j)=0$, $j=1,\dots,m$;
$\mathcal{C}_{\ell}(\eta^*_k)=\delta_{k\ell}$, $k,\ell=1,\dots,n$.  The
$2$-differentials $\eta_k=\frac{-t_k}{\pi}\eta^*_k,\,k=1,\dots,n$ have the
desired pairings with $\frac{\partial}{\partial t_k}$  and the desired expansion.  The final matter is
to note that the pairings of $\{\varphi_j,\psi_j,\eta_k^*\}$  
with $\{\frac{\partial\hat{\nu}}{\partial
s_{\ell}},\frac{\partial\hat{\nu}}{\partial \bar s_{\ell}}\}$ are 
continuous in $(s,t)$. By construction the differential $\hat{\nu}(s)$ is
supported in the complement of the plumbing collars \cite[Lemma 1.1]{Wlspeclim}.  On the support of $\hat{\nu}(s)$ the $2$-differentials are
real analytic in $(s,t)$.  The pairings are continuous and even real analytic. 
\end{proof}

We note general matters for consideration below.  As above, for
$z$ a plumbing collar coordinate for $R_{s,t}$, a quadratic differential
$\psi$  can be factored on the collar as the product of
$\bigl(\frac{dz}{z}\bigr)^2$ and a holomorphic function.  The value $\mathcal{C}(\psi)$
is the constant coefficient of the Laurent expansion of the function.
The surface $R_{s,t}$ is constructed by plumbing $\bigl(R_s\bigr)^*_{c_*^2}$
with the $R_s$-hyperbolic cusp coordinates.  The surface $R_{s,t}$ is the {\em disjoint union}
of $\bigl(R_s\bigr)^*_{c_*}$, $R_s$ with the cuspidal discs $|z_*|,|w_*|<c_*$
removed, and annuli $\{|t|/c_*<|z|<c_*\}$.  

We further note a general matter, the metric-dual map for a Hermitian product space.  For a vector space $V$ 
with Hermitian product $\bigl<\ ,\ \bigr>$, a complex conjugate-linear isomorphism to the dual space $V^*$ 
is defined by the rule $v\in V$ is mapped to the linear functional $\delta(v)=\bigl<\ ,v\bigr>\in V^*$.  
The dual Hermitian product $\bigl<\bigl<\ ,\ \bigr>\bigr>$ for $V^*$ is determined for $v,w\in V$ by $\bigl<\bigl<\delta(v),\delta(w)\bigr>\bigr>=\bigl<w,v\bigr>$.  
For $\{v_j\}$ a basis (over $\mathbb R$ or $\mathbb C$) for $V$  the dual basis for $V^*$ can be expressed in terms of $\delta$.  For the Hermitian pairing matrix $(g_{j\ell})=(\bigl<v_j,v_{\ell}\bigr>)$ and the inverse matrix $(g^{j\ell})$, then $\{\delta(\Sigma_{\ell}g^{j\ell}v_{\ell})\}$ is the corresponding dual basis for $V^*$.  In particular the evaluation of $v_k$ with $\delta(\Sigma_{\ell}g^{j\ell}v_{\ell})$ is $\Sigma_{\ell}\overline{g^{j\ell}}g_{k\ell}=\delta_{kj}$.  
For the dual basis $\{v_j^*\}$ of $V^*$ the pairing matrix has entries $\bigl<\bigl<v_j^*,v_k^*\bigr>\bigr>=\bigl<\bigl<\delta(\Sigma_{\ell}g^{j\ell}v_{\ell}),\delta(\Sigma_m g^{km}v_m)\bigr>\bigr>=\Sigma_{\ell,m}g^{km}g_{m\ell}\overline{g^{j\ell}}=g^{kj}$. The pairing matrix for the dual 
metric is the conjugate inverse of the original pairing matrix.  

We note the particulars for Teichm\"{u}ller space 
and the WP metric.  At a marked Riemann surface $R$ the tangent space of Teichm\"{u}ller space is naturally 
identified with $\mathcal H(R)$, the space of Beltrami differentials harmonic relative to the $R$ hyperbolic 
metric $ds^2$ \cite{Ahsome}.  An element $\mu$ of $\mathcal H(R)$ is given as $\mu=\overline{\varphi}(ds^2)^{-1}$ for $\varphi=\overline{\mu}(ds^2)\in Q(R)$, an integrable holomorphic quadratic differential. The pairing of $\mathcal H(R)$ 
and $Q(R)$ is given by integration.  For 
$\varphi, \psi \in Q(R)$ the integral $\int_R\psi\overline{\varphi}(ds^2)^{-1}$ 
is both the cotangent-tangent pairing of $\psi$ and $\overline{\varphi}(ds^2)^{-1}$ and the Hermitian pairing 
of $\psi$ and $\varphi$.  The mapping $\varphi\rightarrow\overline{\varphi}(ds^2)^{-1}$ is 
the conjugate-linear metric-dual mapping from $Q(R)$ to $\mathcal H(R)$.  For $\{\mu_j\}$ 
a basis for $\mathcal H(R)$, $\{\kappa_j\}$ the dual basis for $Q(R)$, and the pairing matrix $(g_{j\ell})=(\bigl<\kappa_j,\kappa_{\ell}\bigr>_{WP})$, the bases are related by $\overline{\Sigma_{\ell}g^{j\ell}\kappa_{\ell}}=\mu_j(ds^2)$.  

We now apply the two-term expansion of Theorem
\ref{hypexp} for the hyperbolic metric to obtain a two-term expansion for the WP metric. 
The metric expansion will be substituted for the original construction of \cite[Sec. 6]{Msext} to obtain a refinement of the original expansion.  In the following the Eisenstein series 
$E_{k,1}$ and $E_{k,2}$ are for the pair of cusps representing the $k^{{th}}$ node.

\begin{theorem}
\label{wp}
For a noded Riemann surface R with punctures the hyperbolic metric plumbing coordinates for
$R_{s,t}$ provide real analytic coordinates for a local manifold cover
neighborhood for $\Mbar$.  The parameterization is holomorphic in $t$ for $s$
fixed. On the local manifold cover the WP metric is formally Hermitian
satisfying: \begin{enumerate}
\item For $t_k=0$, $k=1,\dots,n$, the restriction of the metric is a smooth K\"{a}hler metric, isometric to the WP product metric for a product of Teichm\"{u}ller spaces.
\item For the tangents $\{\frac{\partial}{\partial s_j},\frac{\partial}{\partial \bar s_j},\frac{\partial}{\partial t_k}\}$ and the quantity $\sigma = \sum\limits^n_{k=1}(\log |t_k|)^{-2}$ then:
$g_{WP}(\frac{\partial}{\partial t_k},\frac{\partial}{\partial t_k})(s,t)=\dfrac{\pi^3}{|t_k|^2(-\log^3 |t_k|)}\,(1\,+\,O(\sigma))$;\\
$g_{WP}(\frac{\partial}{\partial t_k},\frac{\partial}{\partial t_{\ell}})\quad  is\quad  O((|t_kt_{\ell}|\log^3|t_k|\log^3|t_{\ell}|)^{-1})\ for\ k \ne \ell$;\\
and for the tangents $\mathfrak{u}= \frac{\partial}{\partial s_j},\,\frac{\partial}{\partial \bar s_j}$:\\ 
$g_{WP}(\frac{\partial}{\partial t_k},\mathfrak{u})\quad is \quad O((|t_k|(-\log^3|t_k|))^{-1})$. \\
\item For $\mathfrak u= \frac{\partial}{\partial s_j},\,\frac{\partial}{\partial \bar s_j}$, represented at $R_{s,0}$ by $\mu_j$ and $\mathfrak{v}= \frac{\partial}{\partial s_{\ell}}\,,\frac{\partial}{\partial \bar s_{\ell}} $ represented at $R_{s,0}$ by $\mu_{\ell}$ then:\\
$g_{WP}(\mathfrak u,\mathfrak v)(s,t)=\\ g_{WP}(\mathfrak u,\mathfrak v)(s,0)
+  \frac{4\pi^4}{3}\sum\limits_{k=1}^n(\log|t_k|)^{-2}\bigl<\mu_j,\mu_{\ell}(E_{k,1}+E_{k,2})\bigr>_{WP}(s,0) \\
+ O(\sum\limits_{k=1}^n(-\log|t_k|)^{-3})$.
\end{enumerate}
\end{theorem}

\begin{proof}We develop the expansion of the dual metric for the basis provided
in Lemma \ref{coords}. In order to consider integrals we decompose $R_{s,t}$ into $\bigl(R_s\bigr)^*_{c_*}$ and the collars $\{|t|/c_*\le|z|\le c_*\}$.  We compute principal terms of expansions and gather remainder terms.  From Lemma \ref{coords} the variation of the $2$-differentials $\varphi_j$, $\psi_j$,$\eta_k$ is a remainder term.  From Theorem \ref{hypexp} for the principal term for the hyperbolic metric we have $dg_t^2(1+\frac{4\pi^4}{3}\sum(\log|t_k|)^{-2}(E_{k,1}+E_{k,2}))$ on $\bigl(R_s\bigr)_{c_*}^*$ and $ds_t^2$ the plumbing variety metric on the collars.  
On a collar $\{|t|/c_*\le|z|\le c_*\}=\{|t|/c_*\le|w|\le c_*\}$
a quadratic differential is considered as the product of
$\bigl(\frac{dz}{z}\bigr)^2=\bigl(\frac{dw}{w}\bigr)^2$ and a function factor.

{\bf The first term.}  We begin with a model calculation for the hyperbolic metric on the plumbing variety

\begin{gather*}
\int_{\{|t|/c_*<|z|<c_*\}}\bigl|z^{\alpha}\bigl(\frac{dz}{z}\bigr)^2\bigr|^2\,\bigl(ds^2_t\bigr)^{-1}=\\
\frac{2}{\pi}\int^{c_*}_{|t|/c_*}\bigl(\log |t| \sin \frac{\pi \log r}{\log |t|}\bigr)^2r^{2\alpha}\,d\log r \\
\mbox{where for }\alpha=0,\,\mu=\frac{\log r}{\log |t|}\mbox{ and }\epsilon=\frac{\log c_*}{\log |t|},\mbox{ then}\\
=\frac{2}{\pi}(-\log^3|t|)\int^{1-\epsilon}_{\epsilon}\sin^2\pi\mu\,d\mu\,=\,\frac{1}{\pi}(-\log^3|t|)\,+\,O(1),\\
\mbox{ and for }\alpha=1, \mbox{ since }|\sin\mu|\le |\mu|, \mbox{ then }\\
=O(1).
\end{gather*}

We first consider the contribution to the pairing from the
collars.  Consider the contribution for the $\ell^{\mbox{\scriptsize{\,th}}}$
collar.  By construction $\eta_{\ell}$ is the unique quadratic differential
from the dual basis with a nonzero $\mathcal{C}_{\ell}$ evaluation; as above $\eta_{\ell}=-\frac{t_{\ell}}{\pi}\eta_{\ell}^*$.  In
particular $\mathcal{C}_{\ell}(\eta^*_{\ell})=1$ and as above the contribution to the
self pairing for $\eta^*_{\ell}$ for the collar is $\frac{1}{\pi}(-\log^3|t_{\ell}|)+O(1)$.
In general a quadratic differential on a collar is
factored as $\bigl(\frac{dz}{z}\bigr)^2(f_z\,+\,\boldsymbol{c}\,+\,f_w)$ for
$f_z$ holomorphic in $|z|<c,\,f_z(0)=0$; $\boldsymbol{c}$ the
$\mathcal{C}$-evaluation value and $f_w$ holomorphic in $|w|<c,\,f_w(0)=0$.
Furthermore $f_z$, resp. $f_w$, is given as the Cauchy integral of $f$ over
$|z|=c_*$, resp. $|w|=c_*$.  From the Schwarz Lemma $|f_z|\le
c'|z|\max_{|z|=c_*}|f|$ with a corresponding bound for $|f_w|$.  
The $2$-differentials $\varphi_j,\,\psi_j,\,\eta^*_k$ on $|z|=c_*+\epsilon_0$ and $|w|=c_*+\epsilon_0$
depending continuously on $(s,t)$ and so admit decompositions with $f_z$ and $f_w$ uniformly bounded.  The uniform bounds are combined with the majorant bound $|\sin\mu|\le|\mu|$ to bound integrands for the collars.

We combine the above considerations and substitute the grafted metric $dg_t^2$ for the metric principal term to obtain the expansions

\begin{gather*}
\begin{split}
\bigl<\eta^*_k,\,\eta^*_k\bigr>_{WP}\,&=\,\frac{1}{\pi}(-\log^3|t_k|)\,(1\,+\,O(\sum\limits^n_{\ell=1}(\log|t_{\ell}|)^{-2})),\\
\bigl<\eta^*_k,\,\eta^*_{\ell}\bigr>_{WP}\,&=\,O(1)\mbox{  for  }k\ne\ell,\\
\mbox{and for }\mathfrak{a}=\varphi_j,\psi_j;\, \mathfrak{b}=\varphi_{\ell},\psi_{\ell}:\\
\bigl<\mathfrak{a},\,\eta^*_k\bigr>_{WP}\,&=\,O(1).
\end{split}
\end{gather*}
The pairing $\bigl<\mathfrak{a},\mathfrak{b}\bigr>_{WP}$ will be considered in the following paragraphs.

{\bf The second term.} We are ready to determine the second-term in the expansion for $\bigl<\mathfrak{a},\mathfrak{b}\bigr>_{WP}(s,t)$ with $\mathfrak{a}=\varphi_j,\psi_j;\, \mathfrak{b}=\varphi_{\ell},\psi_{\ell}$ and the contribution of the Eisenstein series. 
From Lemma \ref{coords} and Theorem \ref{hypexp} a main step is to give an expansion in $t$ for the integrals of $\mathfrak a\overline{\mathfrak b}(dg_t^2)^{-1}$ and
$\mathfrak a\overline{\mathfrak b}(E^{\dagger}_{k,1}+E^{\dagger}_{k,2})(dg_t^2)^{-1}$ over a collar $\{|t|/c_*\le |z|,|w|\le c_*\}$.  The goal is expansions with remainders at least $O((-\log |t|)^{-3})$ and principal terms with remainders $O(c_*)$ for $c_*$ tending to zero.  

The expansion for $\mathfrak a\overline{\mathfrak b}$ is the first consideration.  The holomorphic quadratic differentials have vanishing $\mathcal C$-evaluations and thus on the collars have factorizations $\mathfrak a=(\frac{dz}{z})^2(\mathfrak a_z+\mathfrak a_w)$ and  $\mathfrak b=(\frac{dz}{z})^2(\mathfrak b_z+\mathfrak b_w)$.  
As above, the summands are bounded as $|\mathfrak a_z|, |\mathfrak b_z|\le |z|c'$ and 
$|\mathfrak a_w|, |\mathfrak b_w|\le |w|c'$ for a suitable constant depending on $\mathfrak a$ and $\mathfrak b$. Now from expansion (\ref{hyp}), the bound $|\sin \mu|\le 1$ and the relation $zw=t$ we have that 
$(|\mathfrak a_z\overline{\mathfrak b_w}|+|\mathfrak a_w \overline{\mathfrak b_z}|)|\frac{dz}{z}|^4(dg_t^2)^{-1}\le 2c'^2|t|(\log|t|)^2|\frac{dz}{z}|^2$.  Furthermore the meldings $E^{\dagger}_{k,*}$ are uniformly bounded in $t$.  It  follows that the integrals of $\mathfrak a_z\overline{\mathfrak b_w}$ and $\mathfrak a_w\overline{\mathfrak b_z}$ over a collar are bounded by $O(|t|(-\log|t|)^3)$ with a similar bound on including the meldings.  
The next matter is the contribution of $\mathfrak a_z\overline{\mathfrak b_z}$ and $\mathfrak a_w\overline{\mathfrak b_w}$; we illustrate by considering the former.  The expansion in $t$ of $\mathfrak a$ and $\mathfrak b$ is 
the first consideration.  The summands $\mathfrak a_z$ and $\mathfrak b_z$ vanish at the origin and are holomorphic in $t$.  The differences $\mathfrak a_z(t)-\mathfrak a_z(0)$ and $\mathfrak b_z(t)-\mathfrak b_z(0)$ are bounded 
by the Schwarz Lemma by $|zt|\le c'$, again for a suitable constant depending on $\mathfrak a$ and $\mathfrak b$.  
It follows that $|\mathfrak a_z(t)\overline{\mathfrak b_z(t)}||\frac{dz}{z}|^4(dg_t^2)^{-1}=
|\mathfrak a_z(0)\overline{\mathfrak b_z(0)}||\frac{dz}{z}|^4(dg_t^2)^{-1}+O(c'^2|t|(\log|t|)^2|dz|^2)$, again with a similar bound on including the meldings. The integral of the remainder term over the collar is $O(|t|(\log|t|)^2)$.  Overall the $z$-term contribution of $\mathfrak a\overline{\mathfrak b}$ has integrand principal term 
$\mathfrak a_z(0)\overline{\mathfrak b_z(0)}|\frac{dz}{z}|^4(dg_t^2)^{-1}$ and remainder term $O(|t|(-\log|t|)^3)$.  

The expansion in $t$ of $(dg_t^2)^{-1}$ is the next consideration.  We proceed by separately considering the subdomains $\{|t|^{1/2}\le |z|\le c_*\}$ and $\{|t|/c_*\le |z| \le |t|^{1/2}\}$. We present an expansion for the first subdomain and an overall bound for the second subdomain.  On the interval $[0,\frac{\pi}{2}]$ we have that $0\le \mu-\frac16\mu^3\le\sin\mu\le\mu$ and thus $0\le \mu^2-\sin^2\mu\le\frac13\mu^4-\frac{1}{36}\mu^6$. We have in general from (\ref{hyp2}) that $ds_0^2\le dg_t^2\le ds_t^2$. It follows for the first subdomain (provided $c_*\le e^{-1}$) that the quantity $(\log|z|)^2-|\frac{dz}{z}|^2(dg_t^2)^{-1}$ is bounded by $O((\log|t|)^{-2}(\log|z|)^6)$.  On the first domain we now have the expansion 
$\mathfrak a_z(0)\overline{\mathfrak b_z(0)}|\frac{dz}{z}|^4(dg_t^2)^{-1}=\mathfrak a_z(0)\overline{\mathfrak b_z(0)}|\frac{dz}{z}|^4(ds_0^2)^{-1}+O(c'^2(\log|t|)^{-2}(\log|z|)^6|dz|^2)$, again with a similar bound on including the meldings. The integral of the remainder term over the first subdomain is $O(c_*(\log|t|)^{-2})$.  We are ready to consider the second subdomain $\{|t|/c_*\le|z|\le |t|^{1/2}\}$ where $|\frac{dz}{z}|^2(dg_t^2)^{-1}$ is bounded by $(\log|t|)^2$ and the integrand $\mathfrak a_z(0)\overline{\mathfrak b_z(0)}|\frac{dz}{z}|^4(dg_t^2)^{-1}$ is bounded by $c'^2(\log|t|)^2|dz|^2$.  The integral over the second subdomain is bounded by $O(|t|(\log|t|)^2)$, again with a similar bound on including the meldings.  

We are ready to gather our considerations and present the expansions for integrals over a collar
\begin{multline*}
\int_{\{|t|/c_*\le|z|\le c_*\}}\mathfrak a\overline{\mathfrak b}(dg_t^2)^{-1}\\
=\int_{\{|t|^{1/2}\le |z|\le c_*\}}\mathfrak a_z(0)\overline{\mathfrak b_z(0)}\bigl|\frac{dz}{z}\bigr|^4(ds_0^2)^{-1}\\+
\int_{\{|t|^{1/2}\le |w|\le c_*\}}\mathfrak a_w(0)\overline{\mathfrak b_w(0)}\bigl|\frac{dw}{w}\bigr|^4(ds_0^2)^{-1} \\
+\ O(c_*(\log|t|)^{-2})\ +\ O(|t|(-\log|t|)^3)
\end{multline*}
and
\begin{multline*}
\int_{\{|t|/c_*\le|z|\le c_*\}}\mathfrak a\overline{\mathfrak b}(E^{\dagger}_{k,1}+E^{\dagger}_{k,2})(dg_t^2)^{-1}\\
=\int_{\{|t|^{1/2}\le |z|\le c_*\}}\mathfrak a_z(0)\overline{\mathfrak b_z(0)}
(E^{\dagger}_{k,1}+E^{\dagger}_{k,2})\bigl|\frac{dz}{z}\bigr|^4(ds_0^2)^{-1}\\+
\int_{\{|t|^{1/2}\le |w|\le c_*\}}\mathfrak a_w(0)\overline{\mathfrak b_w(0)}(E^{\dagger}_{k,1}+E^{\dagger}_{k,2})
\bigl|\frac{dw}{w}\bigr|^4(ds_0^2)^{-1}
+\ O((\log|t|)^{-2}).
\end{multline*}
Observations are in order.  The difference between a melding $E^{\dagger}$ and an Eisenstein series $E$ is bounded by $c(\log|z|)^2$.  The meldings can now be replaced by Eisenstein series with the additional contributions bounded by $O(c_*)$.  The resulting integrals on the right hand side over $\{|t|^{1/2}\le |z|,|w|\le c_*\}$ can now be extended to integrals over the discs $\{|z|,|w|\le c_*\}$ with the additional contributions bounded by $O(|t|(\log|t|)^4)$.  The resulting integrals on the right hand side are now for functions on domains in $R_{s,0}$.

{\bf The dual metric expansion.} We now combine the expansions for the integrals over collars with Theorem \ref{hypexp} to find the expansion for the remaining pairing 
\begin{multline*}
\bigl<\mathfrak a,\mathfrak b\bigr>_{WP}(s,t)=\bigl<\mathfrak a,\mathfrak b\bigr>_{WP}(s,0)\\
-\frac{4\pi^4}{3}\sum\limits^n_{k=1}(\log|t_k|)^{-2}\bigl(\bigl<\mathfrak a,\mathfrak b(E_{k,1}+E_{k,2})\bigr>_{WP}(s,0)+O(c_*)\bigr)\\
+ O(\sum\limits^n_{k=1}(-\log|t_k|)^{-3}).
\end{multline*}
Observations are again in order.
The coefficients of $(\log|t_k|)^{-2}$ are independent of $c_*$ since the Hermitian products $\bigl<\mathfrak a,\mathfrak b\bigr>_{WP}$ are independent of $c_*$ and the remainder term has strictly smaller order of magnitude in $t$.  It follows that the $O(c_*)$ terms vanish.

{\bf The metric expansion.}  The final matter is to relate the expansion for the dual metric to the expansion for the metric.  The first term is an immediate consequence of the following Lemma.  We consider the second term and denote by $V$ the span of the cotangent fields $\{ds_j,\,d\bar s_j,\,ds_{\ell},\,d\bar s_{\ell}\}$.  For a vector space $V$ with basis $\{v_j\}$ and a family of Hermitian pairings $\bigl<\ ,\ \bigr>_{\epsilon}$, let $A_{\epsilon}=(\bigl<v_j,v_{\ell}\bigr>_{\epsilon})$ be the variation of the pairing matrix.  For a variation $A_{\epsilon}=A+\epsilon B$, the inverse matrix satisfies  
$(A+\epsilon B)^{-1}=A^{-1}-\epsilon A^{-1}BA^{-1}+O(\epsilon^2)=A^{-1}-\epsilon A^{-1}B\overline{(A^t)^{-1}}+O(\epsilon^2)$.  As already noted the pairing matrix for the dual metric is given by the conjugate inverse matrix.  We apply the considerations for $\{\kappa_{\ell}\}$ the basis for $Q(R)$ with pairing matrix $g_{j\ell}=(\bigl<\kappa_j,\kappa_{\ell}\bigr>_{WP})$ to find the pairing for the dual basis  $\mu_j(ds_0^2)=\overline{\Sigma_{\ell}g^{j\ell}\kappa_{\ell}}$ for $\{\mu_j\}$ the dual basis for $\mathcal H(R)$. In conclusion for 
$\mathfrak u=\frac{\partial}{\partial s_j},\,\frac{\partial}{\partial \bar s_j}$ represented at $R_{s,0}$ by $\mu_j$, and $\mathfrak v=\frac{\partial}{\partial s_{\ell}},\,\frac{\partial}{\partial \bar s_{\ell}}$ represented at $R_{s,0}$ by $\mu_{\ell}$ we have the expansion 
\begin{multline*}
g_{WP}(\mathfrak u,\mathfrak v)(s,t)\\=g_{WP}(\mathfrak u,\mathfrak v)(s,0)
+ \frac{4\pi^4}{3}\sum\limits_{k=1}^n(\log|t_k|)^{-2}\bigl<\mu_j,\mu_{\ell}(E_{k,1}+E_{k,2})\bigr>_{WP}(s,0)\\
+ O(\sum\limits_{k=1}^n(-\log|t_k|)^{-3}).
\end{multline*}   
\end{proof}

For $\mathcal{A}$ a symmetric $m+n\times m+n$ matrix

\begin{equation*}
\begin{aligned}
\begin{pmatrix}
\lambda_1	&	\dots		&a_{\ell j}	&\dots	\\
\vdots	&	\lambda_k	&\vdots	&\vdots	\\
a_{j\ell}	&	\dots		&\ddots	&\dots	\\
\dots		&	\dots		&\dots	&	B
\end{pmatrix}
&
\begin{aligned}
&\quad\mbox{with }\lambda_k,\ 1\le k \le n; \\ &\quad a _{j\ell},\ 1\le j\le m+n,\ j\ne \ell,\ 1\le \ell\le n \mbox{ and} \\ &\quad B\,=\,(b_{j\ell})\mbox{ a symmetric }m\times m \mbox{ matrix},
\end{aligned}
\end{aligned}
\end{equation*}
we consider the situation that $\lambda_1,\dots,\lambda_n$ are large compared to the $a_{j\ell}$ and $b_{j\ell}$.  We have the following expansion presented in \cite[Prop. 3]{Wlcomp}.

\begin{lemma}
For\ $\det B\ne0$, and $\rho = \sum\limits^n_{k=1}\lambda_k^{-1}$ then:
\[\det \mathcal{A}=\det B\prod\limits^n_{k=1}\lambda_k(1+O(\rho))\]
and $\mathcal{A}^{-1}=(\alpha_{j\ell})$ where: for  $1\le k\le n,\ \alpha_{kk}=\lambda_k^{-1}(1+O(\rho))$; for $1\le j<\ell\le n,\ \alpha_{j\ell}$ is $O((\lambda_j\lambda_{\ell})^{-1})$; for $1\le j\le n<\ell\le m+n,\ \alpha_{j\ell}$ is $O(\lambda_j^{-1})$, and for $1\le j,\,\ell\le m,\ \alpha_{j+n\,\ell+n}=b^{j\ell}(1+O(\rho))$.  
The constants for the $O$-terms are bounded in terms of $m+n$, $\det B^{-1}$ and $\max \{|a_{j\ell}|,|b_{j\ell}|\}$.
\end{lemma}

By way of interest we recall the normal form for the metric
$dg_{WP}^2$; the result is an immediate consequence of the above Theorem \cite[Cor. 4]{Wlcomp}.

\begin{corollary} \label{wpnormal} For the prescribed hyperbolic metric
plumbing coordinates: 
\[
dg^2_{WP}(s,t)\,=\,\bigl(dg^2_{WP}(s,0)\,+\,\pi^3\sum\limits^n_{k=1}(4dr_k^2+r_k^6d\theta^2_k)\bigr)\,(1+O(\|r\|^3))
\] for $r_k=(-\log|t_k|)^{-1/2},\ \theta_k=\arg t_k$ and $r=(r_1,\dots,r_n)$.\
\end{corollary}


\nocite{Obit1,Obit2,OTW,Wlgeom}

\end{document}